\documentclass[12pt]{amsart}

\def \ft{{\mathfrak t}}

\def \R{{\mathbb R}}

\def \Z{{\mathbb Z}}
\def \C{{\mathbb C}}
\def \CP{{\mathbb C}{\mathbb P}}

\def \N{{\mathbb N}}

\newtheorem{theorem}{Theorem}[section]
\newtheorem{lemma}[theorem]{Lemma}
\newtheorem{proposition}[theorem]{Proposition}

\newtheorem{definition}[theorem]{Definition}
\newtheorem{rem}[theorem]{Remark}
\theoremstyle{remark}
\newtheorem{question}[theorem]{Question}
\newtheorem{exmp}[theorem]{Example}

\numberwithin{equation}{section}

\begin{document}

\baselineskip=16pt

\title[Toric degenerations of weight varieties]{Toric degenerations of weight varieties and applications}

\author[P. Foth, Yi Hu]{Philip Foth and Yi Hu}

\address{Philip Foth  \newline
Department of Mathematics, University of Arizona, Tucson, AZ
85721-0089, USA}

\email{foth@math.arizona.edu}

\address{Yi Hu \newline
Department of Mathematics, University of Arizona, Tucson, AZ
85721-0089, USA \newline
 Center for Combinatorics, IPMC, Nankai
University, Tianjin 300071, China}

\email{yhu@math.arizona.edu}



\date{May 01, 2004}

\begin{abstract} We show that a weight variety, which is
a quotient of a flag variety by the maximal torus, admits a flat
degeneration to a toric variety. In particular, we show that the
moduli spaces of spatial polygons degenerate to polarized toric
varieties with the moment polytopes defined by the lengths of their
diagonals. We extend these results to more general
Flaschka-Millson hamiltonians on the quotients of products of
projective spaces. We also study boundary toric divisors and
certain real loci.
\end{abstract}

\maketitle


\section{Introduction}

Let $G$ be a complex connected semisimple Lie group and let $X$ be
a flag variety of $G$, parameterizing parabolic subgroups of a
given type. Several authors studied flat degenerations of $X$ to a
toric variety, starting with the work of Gonciulea and Lakshmibai
\cite{GL}. Caldero \cite{Caldero} used Kashiwara-Lusztig's
canonical bases and their string parameterization in his
construction, which was later extended by Alexeev-Brion
\cite{AlexBrion} to the case of spherical varieties. In this paper
we use their methods to construct flat degenerations of weight
varieties, which are, by definition \cite{Knutson}, quotients of
the flag varieties by the action of the maximal torus. The
resulting polytopes are certain ``slices'' of the string
polytopes.

Interesting examples of weight varieties are the quotients of
complex grassmannians. For example, the quotients of ${\rm
Gr}_{\C}(2, n)$ can be identified \cite{HK} with the moduli spaces
of spatial $n$-gons studied by Klyachko \cite{Klyachko},
Kapovich-Millson \cite{KM}, and many other authors. There are
remarkable integrable systems on these spaces, where the action
variables are given by the lengths of diagonals emanating from a
fixed vertex and the angle variables define the so-called bending
flows, which have a transparent geometric meaning \cite{KM}. As an
application, we show that there exist flat degenerations of the
moduli spaces of polygons to polarized toric varieties, {\it toric
polygon spaces}, whose moment polytopes are defined by the action
variables. These action variables can be computed using the
Gelfand-Tsetlin functions and the Gelfand-MacPherson
correspondence \cite{HK}. Our results can also be generalized for
other grassmannians and the bending flows defined by the
Flaschka-Millson hamiltonians \cite{FlMill}.

We show that there are some real cycles of codimension 2 on the
moduli spaces of polygons which degenerate to the toric
subvarieties corresponding to the facets of the polytopes. We also
extend the work of Kamiyama-Yoshida \cite{Kamiyama} and construct
topological spaces with compact torus action, which conjecturally
are homeomorphic to the central fibers of the flat families.

In the last section of the paper we consider real loci of the
aforementioned spaces and show that they map surjectively to the
moment polytopes and compute the cardinality of the fibers of
these maps.

Let $\overline{M}_{0,n}$ be the moduli space of stable $n$-pointed
rational curves. It can be realized \cite{Kapranov} as the Chow
quotient of the Grassmannian ${\rm Gr}_{\C}(2,n)$. For an
alternative construction using stable polygons, see \cite{Hu99}.
Then $\overline{M}_{0,n}$ can be flatly degenerated to a toric
moduli space whose associated fan is the common refinement of the
fans of all the toric polygon spaces. This and related topics will
appear in a forthcoming paper.

\section{Toric degenerations of weight varieties.}

In this section we use toric degenerations of flag varieties
constructed by Caldero \cite{Caldero} to construct toric
degenerations of weight varieties \cite{Knutson}, which are
defined as GIT quotients of flag varieties by the action of the
maximal torus.

Let $G$ be a connected complex semisimple group, $B$ a Borel
subgroup, $U$ its unipotent radical, and $H$ a Cartan subgroup
such that $B=HU$. Let also $\Phi=\Phi(G, H)$ be the system of
roots, $\Phi^+=\Phi^+(B, H)$ the subset of positive roots, and $\{
\alpha_1, ...,\alpha_r\}$ the basis of simple roots, where $r$ is
the rank of $G$. Let $\Lambda$ be the weight lattice of $G$ and
$\Lambda^+$ the subset of dominant weights. For
$\lambda\in\Lambda^+$ we denote by $V(\lambda)$ the irreducible
$G$-module with highest weight $\lambda$. Let $P_\lambda\supset B$
be the parabolic subgroup of $G$ which stabilizes a highest weight
vector in $V(\lambda)$. Also denote by $L_\lambda = G
\times_{P_\lambda} \C$ the $G$-linearized line bundle on
$X_\lambda:=G/P_\lambda$ corresponding to the character $\lambda$
extended to $P_\lambda$.

Let $W$ be the Weyl group and $w_0\in W$ the longest element of
length $\ell$. Choose a reduced decomposition
$$
\underline{w_0}=s_{i_1}s_{i_2}\cdots s_{i_\ell}
$$
into a product of simple reflections. Recall that there is a
rational convex polyhedral cone ${\mathcal
C}_{\underline{w_0}}\subset \Lambda_\R\times \R^\ell$ such that
the integral points in the {\it string polytope} $Q(\lambda)$
correspond to the elements of the (dual) canonical basis (for the
chosen string parameterization) $(b_{\lambda,\phi})$ for any fixed
$\lambda\in\Lambda^+$. For a fixed $\lambda$, we may identify
$Q(\lambda)$ with its image in $\R^\ell$ via the projection
$\Lambda_\R \times \R^\ell \to \R^\ell$.  If we denote the string
parameterization by
$$
b_{\lambda,\phi}\mapsto (\lambda, t_1, ..., t_\ell)\in
\Lambda^+\times\N^\ell,
$$
then the projection
$$
\pi_\lambda:\ \R^\ell\to \Lambda_\R, \ \ (t_1, ..., t_\ell)\mapsto
-\lambda+t_1\alpha_{i_1}+\cdots +t_\ell\alpha_{i_\ell}
$$
maps the string polytope $Q(\lambda)$ onto the convex hull of the
Weyl group orbit of the dual weight $\lambda^*=-w_0\lambda$:
$$
\pi_\lambda(Q(\lambda))={\rm Conv}(W.\lambda^*)=-{\rm
Conv}(W.\lambda):=\Delta(\lambda)
$$

Caldero \cite{Caldero} (see also \cite{AlexBrion}) has constructed
a flat deformation of the polarized flag variety $(X_\lambda,
L_\lambda)$ to a polarized toric variety $(X_{\lambda; 0},
L_{\lambda; 0})$ such that the corresponding moment polytope is
exactly the string polytope $Q(\lambda)$. His construction is
based on the key multiplicative property of the (dual) canonical
basis:
$$
b_{\lambda_1, \phi_1}b_{\lambda_2, \phi_2}=b_{\lambda_1+\lambda_2,
\phi_1+\phi_2}+ \sum_{\phi < \phi_1+\phi_2} {\rm
coeff}.b_{\lambda_1+\lambda_2, \phi},
$$
which allows to put a filtration on the homogeneous coordinate
ring which is compatible with grading and the torus action such
that the associated graded ring is the algebra of
the semigroup of integral points in a rational convex polyhedral
cone, see \cite{Caldero}, \cite{AlexBrion} for details.

Let $\Phi_\lambda: X_{\lambda; 0} \rightarrow Q(\lambda) \subset
\R^\ell$ be the moment map for ${\mathbb T}$-action, where 
${\mathbb T}$ is the compact part of $(\C^*)^\ell$. Then the
composition $\pi_\lambda \circ \Phi_\lambda: X_{\lambda; 0} \to
\Delta_\lambda \subset \Lambda_\R$ is the moment map for the torus
$H$.

Recall that a {\it weight variety} $M_{\lambda; \mu}$ is a GIT
quotient of $X$ by the action of $H$ associated with an integral
point $\mu \in \Delta(\lambda)$. More precisely, let $\mu \in
\Delta(\lambda)$ be a character of $H$, and $L_\lambda (-\mu)$ be
the linearized line bundle $L_\lambda$ twisted by $-\mu$. Then
$M_{\lambda; \mu} = X_\lambda^{ss}(L_\lambda (-\mu))/\!/G$. In
symplectic terms, the $H$- moment map for the polarized variety
$(X_\lambda,L_\lambda (-\mu))$ is $\pi_\lambda \circ \Phi_\lambda
- \mu$, and $M_{\lambda; \mu}$ is identified with the reduction of
$X_\lambda$ at the level $\mu$. The line bundle $L_\lambda (-\mu)$
descends to the quotient $M_{\lambda; \mu}$. Let us denote the
descended line bundle by $L_{\lambda;\mu}$.

\begin{theorem} There exists a flat degeneration of the weight variety
$M_{\lambda; \mu}$ to a projective toric variety $N_{\lambda;
\mu}$ with a polarization $L_{\lambda;\mu;0}$ such that the
deformation carries $L_{\lambda;\mu}^n$ to $L_{\lambda;\mu;0}^n$
for sufficiently large positive integer $n$. In particular, the
moment polytope of the polarized toric variety $(N_{\lambda;
\mu},L_{\lambda;\mu;0})$ is $Q_{\lambda;
\mu}=\pi_{\lambda}^{-1}(\mu) \cap Q(\lambda)$.
\label{th:flat_weight}
\end{theorem}

\proof  By the Borel-Weil Theorem, we have the following
identification of the section ring
$$R_\lambda = \bigoplus_{n=0}^\infty H^0(X, L_\lambda^n) = \bigoplus_{n=0}^\infty
V(n\lambda^*).$$ $H$ acts on $R_\lambda$ via the induced action
twisted by $-\mu$. Let $R_\lambda^{H(-\mu)}$ denote the set of
invariants of this twisted $H$-action. Using the canonical bases,
we have
$$R_\lambda  = \bigoplus_{n=0}^\infty \bigoplus_\phi \C .b_{n
\lambda, \phi}.$$ Since $b_{n \lambda, \phi}$ is an eigenvector of
eigenvalue $- n \lambda + \sum_k t_{i_k} \alpha_{i_k}$ with
respect to the left $H$-action,   we see that after twisting the
$H$ action on $\C .b_{n \lambda, \phi}$ by $-n\mu$, the eigenvalue
becomes zero if $\sum_k t_{i_k} \alpha_{i_k} = n (\lambda + \mu)$.
Hence,
$$R_\lambda^{H (-\mu)} =\bigoplus_{n=0}^\infty \bigoplus_{\sum t_{i_k} \alpha_{i_k}  = n
(\lambda + \mu)} \C .b_{\lambda, \phi}.$$  Observe that the
filtration on $R_\lambda$ induces a filtration on $R_\lambda^{H
(-\mu)}$. Let $\hbox{Gr}(R_\lambda^{H (-\mu)})$ be the associated
graded algebra of  $R_\lambda^{H (-\mu)}$. Then the weight variety
$M_{\lambda; \mu} = \hbox{Proj} (R_\lambda^H (-\mu))$ degenerates
flatly to $N_{\lambda; \mu} =\hbox{Proj} (\hbox{Gr} R_\lambda^{H
(-\mu)})$, and the degeneration carries the set of ample line
bundles $L_{\lambda;\mu}^n$ to the set of  ample line bundles
$L_{\lambda;\mu;0}^n$ on $N_{\lambda, \mu}$ for sufficiently large
$n$ (\cite{AlexBrion}). Since in general, the moment map of a
polarization $L$ is the moment map of $L^n$ divided by $n$,  it
follows from \cite{KSZ} that the moment polytope of the polarized
toric variety $(N_{\lambda; \mu}, L_{\lambda;\mu;0})$ is
$\pi_\lambda^{-1}(\mu) \cap Q(\lambda)$.
\endproof

\begin{definition}
We may call the central toric fiber $N_{\lambda, \mu}$ the toric
weight variety of type $(\lambda, \mu)$, and call the polytope
$Q_{\lambda; \mu} = \pi_\lambda^{-1}(\mu) \cap Q(\lambda)$ the
weight polytope of type $(\lambda, \mu)$. Note that the polytope
$Q_{\lambda; \mu}$ is rational but needs not to be integral in
general.
\end{definition}

\begin{rem}{\em
Unlike the isomorphic type of the quotient $M_{\lambda; \mu}$
which only depends on the GIT chamber of $\mu$, the degeneration
depends on the choice of individual $\mu$, in particular, the
topological type of the central fiber $N_{\lambda; \mu}$ can
change even when $\mu$ varies within its GIT chamber. See Example
\ref{M11111}.}
\end{rem}

\begin{question}
Is there a chamber structure on $\Delta_\lambda$ so that
$N_{\lambda; \mu}$,  as (unpolarized) toric variety, remain the
same within a chamber?
\end{question}

\begin{rem}{\em As noted in \cite[Section 5]{AlexBrion},
the Ehrhart polynomial of the polytope $Q_{\lambda;\mu}$ is given
by
$$ n\mapsto {\rm dim}V(n\lambda^*)_{n\mu}, $$ the multiplicity
of the weight $n\mu$ for the right $H$-action in the irreducible
$G$-module with highest weight $n\lambda^*$. See also discussion
on Pieri's formula in \cite{FlMill}.}
\end{rem}

\section{Degeneration of polygon spaces}

Let us recall the generalities about the polygon spaces. We will
start with spatial polygons in Eucledian 3-space. Given an
$n$-tuple of positive real numbers ${\bf r}=(r_1, ..., r_n)$ we
consider the moduli space $M_{\bf r}$ of $n$-gons with consecutive
side lengths $r_1, ..., r_n$. This space can be viewed as the
symplectic quotient at the zero level of the product of
two-dimensional spheres $ (S^2)^n$ by the diagonal action of the
group ${\rm SO}(3)$. The symplectic structure on the $j$-th
multiple is taken to be $r_j$ times the standard unit sphere
volume form. On the algebraic geometry side, the space $M_{\bf r}$
can be realized as the GIT quotient of the $n$-fold product of the
projective line $\CP^1$ by the diagonal action of the group ${\rm
SL}(2, \C)$. The choice of linearization is given by $r_i$'s
(which we assume to be integers here). These spaces were studied
by Klyachko \cite{Klyachko}, Kapovich-Millson \cite{KM},
Hausmann-Knutson \cite{HK}, and many others.

Let us now recall the Gelfand-MacPherson correspondence
\cite{GelMac} in the form of Theorem 2.4.7 of \cite{Kapranov}. It
states that for any $n$-tuple of positive integers ${\bf r}=(r_1,
..., r_n)$, there is an isomorphism of GIT quotients
$$
({\rm Gr}(k, n)/\!/H_1)_{{\mathcal O}(1), {\bf r}}\simeq
((\CP^{k-1})^n/\!/{\rm SL}(k, \C))_{{\mathcal O}({\bf r}), \zeta},
$$
where the line bundle ${\mathcal O}(1)$ on the grassmannian is
linearized to correspond to the action of $H_1=\{ (t_1, ...,
t_n)\in (\C^*)^n: \prod t_i=1\}$ on $\C^n$ given by
$$
(t_1, ..., t_n)\mapsto {\rm diag}(t^{\bf r}t_1, ..., t^{\bf
r}t_n),
$$
and $t^{\bf r}=t_1^{r_1}\cdots t_n^{r_n}$ is the character of
$H_1$ corresponding to ${\bf r}$. On the other side, the line
bundle ${\mathcal O}({\bf r})$ on $(\CP^{k-1})^n$ is the tensor
product of pull-backs of line bundles ${\mathcal O}(r_i)$ from the
$i$-th multiple. This bundle has exactly one ${\rm SL}(k,
\C)$-linearization, which is denoted by $\zeta$.

Now, let us apply this correspondence to $M_{\bf r}$ using the
language of hermitian matrices as was first done in \cite{HK}. Let
us start with the complex vector space $\C^{2n}$ with coordinates
$(z_1, ..., z_n, w_1, ...,w_n)$. For convenience, we arrange
$z_i$'s and $w_j$'s in a matrix form
$$
A=\left(
\begin{array}{cccc}
z_1 & z_2 & \cdots & z_n \\
w_1 & w_2 & \cdots & w_n
\end{array}
\right)
$$
If we let $P=\sum_{i=1}^n r_i/2$ stand for half the perimeter,
then the moduli space $M_{\bf r}$ can be realized as the quotient
of the grassmannian ${\rm Gr}(2, n)$ by the action of the maximal
torus. The proof goes via the so-called reduction in stages. The
grassmannian  ${\rm Gr}(2, n)$ can be realized as the symplectic
quotient of the space $\C^{2n}$ with the standard Darboux
symplectic form by the action of the group ${\rm U}(2)$ by left
multiplication on the matrix $A$. The corresponding moment map
takes values in $2\times 2$ hermitian matrices and is given by
$$
A{\overline{A}}^t=\left(
\begin{array}{cc}
\sum_{i=1}^n |z_i|^2 & \sum_{i=1}^n z_i{\bar w}_i \\ {} & {} \\
\sum_{i=1}^n w_i{\bar z}_i & \sum_{i=1}^n |w_i|^2
\end{array}
\right).
$$
And we reduce at the level ${\rm diag}(P, P)$. Alternatively, we
can reduce the space $\C^{2n}$ by the action of $(S^1)^n$ on the
columns of the matrix $A$ at the level $(r_1, ..., r_n)$ to get
the n-fold product of two-dimensional spheres with the
aforementioned symplectic structure.

Using the $2\times 2$ hermitian matrices, we see that the moduli
space of polygons can be identified with the space of solutions of
matrix equation
$$
A_1+\cdots + A_n=\left( \begin{array}{cc} P & 0 \\ 0 & P
\end{array}\right),
$$
where each matrix $A_j$ has the prescribed spectrum $(0, r_j)$ and
in the language of the above correspondence:
$$
A_j=\left( \begin{array}{cc} |z_j|^2 & z_j{\bar w}_j \\ {} & {} \\
{\bar z}_jw_j & |w_j|^2 \end{array}\right).
$$

Let $d_j$ stand for the length of the diagonal of the polygon
connecting vertices $1$ and $j$. There is a remarkable integrable
system \cite{KM} on $M_{\bf r}$ whose action variables are these
functions $d_j$ for $3\le j\le n-1$. The corresponding flows are
the so-called {\it bending} flows and have a simple geometric
description by rotating the last $n-j$ spheres about the axis
defined by the diagonal $d_j$. This integrable system was shown
\cite{HK} to be the reduction of the so-called Gelfand-Tsetlin
integrable system on ${\rm Gr}(2, n)$ \cite{GS}, which we recall
now briefly.

Consider the $n\times n$ hermitian matrix $A^*A$ of rank $2$. The
two non-trivial eigenvalues of the $k\times k$ upper-left
submatrix denote by $a_j$ and $b_j$, $a_j\le b_j$. (For
completeness, we let $a_1=0$.) The value $b_{n-1}=P$ is
predetermined by the fact that the two non-zero eigenvalues of
$A^*A$ are both equal to $P$. Thus we have $2(n-2)$ independent
commuting hamiltonians on ${\rm Gr}(2,n)$, which are also called
the Gelfand-Tsetlin variables and commonly written in a triangular
form:
$$
\begin{array}{ccccccccccc}
0 & \ & 0 & \ & \cdots & \ & 0 & \ & P & \ & P \\
\ & 0 & \ & 0 & \cdots & 0 & \ & a_{n-1} & \ & P & \ \\
\ & \ & \ & \cdots & \ & \cdots & \ & \cdots & \ & \ & \ \\
\ & \ & \ & 0 & \ & a_3 & \ & b_3 & \ & \ & \ \\
\ & \ & \ & \ & a_2 & \ & b_2 & \ &  \ & \ & \ \\
\ & \ & \ & \ & \ & b_1 & \ & \ & \ & \ & \
\end{array}
$$
which also satisfy the interlacing inequalities.

The reduction of this system corresponds to choosing the diagonal
values of the matrix $A^*A$ to ensure the compatibility with
chosen side length, which amounts to fixing the sums in rows of
the Gelfand-Tsetlin diagram: $b_1=r_1$, $a_2+b_2=r_1+r_2$, ...,
$a_j+b_j=\sum_{i=1}^jr_i$, ... , $a_{n-1}=P-r_n$. Note that the
lengths of the diagonals can also easily be expressed in terms of
the variables $(a_i, b_i)$, namely: $$d_{j+1}=b_j-a_j.$$ These
variable diagonal lengths $d_j$'s, $3\le j\le n-1$, which are the
action variables define a polytope $\Pi_{\bf r}$ in $\R^{n-3}$.
This polytope is cut by $3(n-2)$ hyperplanes corresponding to the
triangle inequalities, some of which are redundant, for the
triangulation of the polygon by all these diagonals emanating from
the first vertex. The $j$-th triangle has the side lengths $d_j,
d_{j+1}$, and $r_{j+1}$ so it contributes three inequalities:
$$d_j+d_{j+1}\le r_{j+1}, d_j+r_{j+1}\le d_{j+1}, d_{j+1}+r_{j+1}\le
d_j. \;\;\;\;\;\;\;\; (*)$$

\begin{definition}
The weight polytope $\Pi_{\bf r}$ is a polytope in $\R^{n-3}$
defined by all the triangle inequalities of ${\rm (*)}$.
\end{definition}

Kogan-Miller \cite{KogMil} and Alexeev-Brion \cite{AlexBrion} have
explicitly constructed degenerations of complex flag varieties $X$
of type $A_n$ to polarized toric varieties which have the
Gelfand-Tsetlin polytopes (denoted by GT$(\lambda)$, $\lambda\in
\Lambda^+$) as their moment polytopes. In the case of
Alexeev-Brion degeneration, they used the simplest reduced
decomposition of the longest Weyl group element in $W=S_{n+1}$
$$
\underline{w}_0^{\rm std}=s_1s_2s_1s_3s_2s_1\cdots
s_ns_{n-1}\cdots s_1.
$$
As before, we denote by $\pi_\lambda: {\rm GT}(\lambda)\to
\Delta(\lambda)\subset\Lambda_\R$ the projection of the
Gelfand-Tsetlin polytope onto the moment polytope of the flag
variety $X$ with respect to the hamiltonian ${\rm SU}(n)$-action
and the invariant symplectic form $\omega$ on $X=G/P_\lambda$ such
that $[\omega]=c_1(L_\lambda)$. From Theorem \ref{th:flat_weight}
we immediately conclude

\begin{theorem}
\label{thm:MtoN} There exists a flat degeneration of the moduli
space of polygons $M_{\bf r}$ with integral side length ${\bf
r}=(r_1, ..., r_n)$ together with a polarizing line bundle $L_r$
on it to a polarized toric variety $N_{\bf r}$ whose moment
polytope is the weight polytope $\Pi_{\bf r}$.
\end{theorem}


The toric variety $N_{\bf r}$, even its topological type,  depends
on the choice of ${\bf r}$. In fact, the combinatorial type of the
weight polytope $\Pi_{\bf r}$ may also change when ${\bf r}$
changes, as shown in the following  examples.

\begin{exmp}
\label{M11111}
 Consider the case of pentagons. Already here we will
see how the degeneration depends on the choice of the linearized
line bundle, which in terms of the polygons corresponds to
choosing not only side length, but also the order in which they
appear. Let us consider the following three quintuples of
integers: ${\bf r}^1=(3,3,3,3,3)$, ${\bf r}^2=(3,3,3,3,4)$, and
${\bf r}^3=(3,4,3,4,3)$. For all these choices, the moduli spaces
$M_{{\bf r}^1}$, $M_{{\bf r}^2}$, and $M_{{\bf r}^3}$ are
biholomorphically isomorphic to the same space
${\overline{M}}_{0,5}$, the canonically compactified moduli space
of projective lines marked at 5 distinct points. This surface is
also isomorphic to the Del Pezzo surface of degree 5 (blow-up of
$\CP^2$ at 4 generic points). The corresponding weight polytops
are  2-dimensional polygons. The coordinates $(x,y)$ correspond to
the lengths of the first and the second diagonals respectively.

In the case of ${\bf r}^1$ we have a pentagon with vertices at the
points $(0,3)$, $(3,0)$, $(6,3)$, $(6,6)$, and $(3,6)$. This
corresponds to a toric surface with two isolated singular points
(local quotient singularities of $\C^2$ by the action of $\Z/2$).

In the case of ${\bf r}^2$ we have a hexagon with vertices at the
points $(0,3)$, $(2,1)$, $(4,1)$, $(6,3)$, $(6,7)$, and $(4,7)$.
This corresponds to a toric surface with one isolated singular
point (over $(3,0)$).

In the case of ${\bf r}^3$ we get a heptagon with vertices at the
points $(1,4)$, $(1,2)$, $(2,1)$, $(4,1)$, $(7,4)$, $(7,7)$, and
$(4,7)$. The corresponding toric variety is smooth (and
diffeomorphic to $\overline{M}_{0,n}$). This toric variety is a
blow up of each of the previous ones at one or two singular
points. Notice that we would get a different polytope with five
vertices should the same side lengths be arranged as
$(3,3,3,4,4)$.
\end{exmp}

\begin{exmp} In the equilateral case when $n=4,5,6$, or $7$ and ${\bf r}=(1,1,...,1)$, the
polytopes $\Pi_{\bf r}$ were considered in \cite[Section
4]{Kamiyama} and their virtual Poincar\'e polynomials were
computed.
\end{exmp}

\begin{exmp} Now consider an example when ${\bf r}=(2,2,2,..., 2, 2n-3)$,
when there is one ``very long'' side. In this case ${\bf r}$ lies
in a {\it favorable} chamber \cite[Section 3]{Hu99}. One can see
that the polytope $\Pi_{\bf r}$ is defined by the following set of
inequalities:
$$
3\le d_3\le 4, \ \ 5\le d_4\le d_3+2, \ \ ..., \ \ 2i-3\le d_i\le
d_{i-1}+2, \ {\rm for} \ 4\le i\le n-1.
$$
This polytope is actually a simplex, and the corresponding
polarized toric variety is $(\CP^{n-3}$, ${\mathcal O}(1))$. This
is not surprising, because for any element ${\bf r}$ from the
interior of a favorable chamber, the moduli space $M_{\bf r}$ is
itself isomorphic to $\CP^{n-3}$, and thus our flat family is
trivial in this case.
\end{exmp}

\begin{exmp} Let us consider the case when we have three ``very long'' sides.
More precisely, let $r_1=r_2=\cdots=r_{n-3}=1$ and
$r_{n-2}=r_{n-1}=r_n=n$. So the perimeter is equal to $4n-3$ and
no two long sides can ever be parallel. We will show that the
moduli space $M_{\bf r}$ is isomorphic to $(\CP^1)^{n-3}$ in a
natural way. Let us treat the moduli space $M_{\bf r}$ as the
moduli space of weighted configurations of points on $\CP^1$.
Denote these points by $z_1$, ..., $z_n$ respectively, where we
coordinatize, as usual, by $\C$ the complement of a point (the
``North pole'' or $\infty$) in $\CP^1$. Consider the cross-ratios
$$
w_i=\frac{z_{n-2}-z_{n}}{z_{n}-z_{n-1}} \cdot
\frac{z_i-z_{n-1}}{z_n-z_i}, \ \ {\rm where} \ \ 1\le i\le n-3.
$$
Since the points $z_{n-2}$, $z_{n-1}$ and $z_n$ never collide on
the subset of semi-stable configurations, these cross-ratios
define a global map from the moduli space $M_{\bf r}$ to
$(\CP^1)^{n-3}$, which is an isomorphism. It is important to
notice that the bending flows are {\bf not} the standard circle
actions on this space.

This example also shows that the toric variety $N_{\bf r}$ and the
polytope $\Pi_r$ depends on the order of side-lengths. If we had
switched the first and the $(n-2)$-nd sides, so that ${\bf r}=(n,
1, 1, ..., 1, n, n)$, then none of the diagonals emanating from
the first vertex would ever vanish, because it would not be
allowed by the triangle inequalities, and the central fiber would
be isomorphic to $(\CP^1)^{n-3}$. But if we stick with the
original choice of order on the side-lengths, then the resulting
toric variety will be singular, in general. This can be seen
already for $n=5$, when the resulting toric variety has an
isolated quotient singularity.
\end{exmp}

In fact, since all the inequalities on the lengths of the
diagonals emanating from a given vertex come from the triangles
which these diagonals break the polygon into, we immediately have
the following result:

\begin{proposition}
\label{facets} The weight polytope $\Pi_{\bf r}$ is the compact
intersection of the half spaces defined by affine hyperplanes of
the following forms: $x_1 = |r_1\pm r_2|, x_{n-3}=|r_n\pm
r_{n-1}|$ or $|x_i \pm x_{i+1}|=r_{i+2}$. In particular, let
$\{e_i\}$ be the coordinate basis vectors, then any facet is
perpendicular to one of the following: $e_1, e_{n-3}$, i.e., the
first or the last coordinate line; $e_i \pm e_{i+1}$, i.e., the
diagonal or anti-diagonal line in the two dimensional coordinate
plane spanned by $\{e_i, e_{i+1}\}$.
\end{proposition}

In the above Proposition $x_k$ corresponds to $d_{k+2}$.

Let $U_{\bf r}$ be the (non-moduli) space of polygons with side
length vector ${\bf r}$ such that the first vertex is at the
origin. Note that $M_{\bf r} =  U_{\bf r}/{\rm SO}(3)$. We construct a
new quotient space by defining a new equivalence relation on
$U_{\bf r}$. For each $3 \le i \le n-1$, let $D_i$ be the set of
polygons with $d_i =0$. Let $D^\circ$ be the set of polygons where
none of the diagonals vanish. Then $$U_{\bf r} = D^\circ \cup_I
D_I,$$ where $I \subset \{3, \ldots, n-1\}$ and $D_I$ is the set
of polygons with vanishing diagonals $d_i, i \in I$ but the rest
of the diagonals do not vanish. For any $I$, every polygon in
$D_I$ can be decomposed as polygons $P_1, \ldots, P_\ell$ mutually
joining at a point, where $\ell=|I|+1$. In this case, set $$S_I =
D_I/({\rm SO}(3))^\ell$$ where $({\rm SO}(3))^\ell$ acts on each factor
component-wise. Let $S^\circ = D^\circ /{\rm SO}(3)$. We then obtain a
new topological space
$$V_{\bf r} = S^\circ \cup_I S_I.$$
There is an obvious continuous collapsing map
$$f: M_{\bf r} \rightarrow V_{\bf r}.$$
There is also a continuous map $$\Phi_{\bf r} : V_{\bf r}
\rightarrow \Pi_{\bf r}.$$ The map $f$ passes the action-angle
coordinates on $M_{\bf r}$ to $V_{\bf r}$.

\medskip\noindent{\bf Conjecture.} The space $V_{\bf r}$ is homeomorphic
to $N_{\bf r}$. \medskip

We were able to verify this conjecture in some particular cases,
however the rigorous proof evades us, because the explicit nature
of the deformation family is evasive.

\begin{rem}{\em
Our construction of the space $V_{\bf r}$ follows Kamiyama-Yoshida
\cite{Kamiyama}, where a similar space a continuous map to the
polytope $\Pi_{\bf r}$ was constructed in the equal side length
case. Despite the fact that they called it a ``symplectic toric
manifold'', their space only had a structure of a topological
space stratified into a finite union of orbifolds.}
\end{rem}

\subsection{Boundary toric divisors}

Corresponding to each facet, there is a toric subvariety of
$N_{\bf r}$ of codimension 1. Corresponding to the set of all
facets perpendicular to a given vector is the union of
subvarieties whose isotropy subgroup is generated by a prime
vector in that direction. By the above proposition, the one
dimensional isotropy subgroups are listed as below:
\begin{enumerate}
\item $T_1= \{t \in (\C^*)^{n-3} | t_2 = \ldots = t_{n-3} = 1\}$.
\item $T_{n-3} = \{t \in (\C^*)^{n-3} | t_1 = \ldots = t_{n-4} =
1\}$. \item $T_i^+ = \{t \in (\C^*)^{n-3} | t_i = t_{i+1}, t_k = 1
\; \hbox{for} \; k \ne i, i+1 \}$. \item $T_i^- =\{t \in
(\C^*)^{n-3} | t_{i+1} = t_i^{-1}, t_k = 1 \; \hbox{for}\; k \ne
i, i+1\}$. .
\end{enumerate}

Corresponding to the list, we have toric divisors (in parentheses
the equations of the corresponding facets of $\Pi_{\bf r}$ are
given):
\begin{itemize}
\item (1) $N^1_2$ ($x_1=r_1+r_2$) and $N^2_2$ ($x_1=r_2-r_1$) or $N^3_2$ ($x_1=r_1-r_2$)
\item (2) $N^3_{n-1}$ ($x_{n-3}=r_{n-1}+r_n$) and $N^1_{n-1}$ ($x_{n-3}=r_n-r_{n-1}$)
or $N^2_{n-1}$ ($x_{n-3}=r_{n-1}-r_n$)
\item (3) $N^1_i$ ($x_{i-1}-x_{i-2}=r_i$) and $N^3_i$ ($x_{i-2}-x_{i-1}=r_i$)
\item (4) $N^2_i$ ($x_{i-1}+x_{i-2}=r_i$)
\end{itemize}
These are all the toric divisors, some of which are reducible. In
the next section, we will show that they are the deformation
images of certain subpaces of $M_{\bf r}$.

The Chow ring $A^*(M_{\bf r})$ is generated, as a ring, by degree
$1$ cycles corresponding to divisors $Z_{ij}$, $i< j$ defined by
the condition that the $i$-th side of the polygon points in the
same direction as the $j$-th side. We call such sides ``positively
parallel''  (or strongly parallel), which is a stronger condition
than just being parallel. This statement can be verified by
considering the cycles $D^S$ defined in \cite{Keel} and which
generate $A^*({\overline{M}}_{0,n})$ together with the proper
surjective birational morphism ${\overline{M}}_{0,n}\to M_{\bf
r}$, and its properties \cite{FothCMB}. It is then well-known that
the map $A^*(M_{\bf r})\to H^{2*}(M_{\bf r}, \Z)$ is an
isomorphism. There are some natural linear relations among the
cycles $Z_{ij}$, which again can be deduced from the linear
relations on the cycles $D^S$ studied in \cite{Keel} and
\cite{KonMan}. However, these cycles are not stable under the
bending flows, and thus will not degenerate to toric subvarieties.

Instead, let us consider subsets of the space $M_{\bf r}$ which
are defined by the condition that the $i$-th diagonal is parallel
to the $(i+1)$-st diagonal, and thus to the $i$-th side  (we also
think that a zero vector is parallel to any other vector). Let us
introduce a notation for these subsets.
\begin{itemize}
\item denote by $Y^1_i$ the subset where $d_i+r_i=d_{i+1}$
\item denote by $Y^2_i$ the subset where $d_i+d_{i+1}=r_i$
\item denote by $Y^3_i$ the subset where $d_i-d_{i+1}=r_i$
\item denote by $Y_i=Y^1_i\cup Y^2_i\cup Y^3_i$ their union
\end{itemize}

Here we have $2\le i\le n-1$ as we think of the first and last
sides as $2$-nd and $n$-th diagonals respectively. The number of
such we get this way is bounded from above by $3n-8$. However,
depending on the choice of ${\bf r}$ some of these subsets might
be empty.

One can easily see that $Y^1_i$ and $Y^3_i$ always have empty
intersection. However, if $d_i$ can attain zero, then $Y^1_i$ and
$Y^2_i$ do intersect and when the side lengths are generic, and
$3<i<n-1$, this intersection is homeomorphic to an $S^1$-bundle
over an $S^2$-bundle over the space $M_{\bf r^{(1)}}\times M_{\bf
r^{(2)}}$, where ${\bf r^{(1)}}=(r_1, ..., r_{i-1})$ and ${\bf
r^{(2)}}=(r_{i}, ..., r_n)$.

Similarly, if $d_{i+1}$ can attain the value of zero, then then
$Y^2_i$ and $Y^3_i$ do intersect and when side lengths are
generic, and $2<i<n-2$, this intersection is homeomorphic to an
$S^1$-bundle over an $S^2$-bundle over the space $M_{\bf
r^{(1)}}\times M_{\bf r^{(2)}}$, where ${\bf r^{(1)}}=(r_1, ...,
r_{i})$ and ${\bf r^{(2)}}=(r_{i+1}, ..., r_n)$.

Note that when $r_1=r_2$ then the subsets $Y^2_{2}$ and $Y^3_2$
coincide and are equivalently defined by $d_3=0$. Same goes for
the case $r_{n-1}=r_n$, subsets $Y^1_{n-1}$ and $Y^2_{n-1}$
coincide and are defined by $d_{n-1}=0$. Unless stated otherwise,
we would like to exclude these from consideration in the statement
of the following result.

We denote by $M_{\bf r}^0$ the ``toric'' open submanifold of
$M_{\bf r}$, where none of the diagonals vanishes. We say that a
subset $Q\subset M_{\bf r}$ {\it degenerates to} $Q_0\subset
N_{\bf r}$ if $Q\cap M_{\bf r}^0$ (resp $Q_0$) is preserved by the
bending flows (resp. compact torus action), and they have the same
image in $\Pi_{\bf r}$.

\begin{lemma}
\begin{enumerate}
\item If $M_{\bf r }$ is smooth, then $Y_{i}\cap M_{\bf r}^0$, for
$2\le i\le n-1$, is a smooth closed submanifold of $M_{\bf r}^0$
of real codimension $2$.
\item Each $Y_i$ degenerates to the union of toric
divisors $N^1_{i}$, $N^2_i$, and $N^3_i$ of $N_{\bf r}$. All toric
divisors  of $N_{\bf r}$ arise this way.
\end{enumerate}
\end{lemma}
\proof (1) The intersections $$Y_{i} \cap M_{\bf r}^0$$ are
clearly non-empty. Each connected components of these
intersections is also a connected component of the fixed point set
of some one-dimensional torus (see the remark after Proposition
\ref{facets}), so it is smooth.

(2) is clear.
\endproof

\begin{rem}{\em The intersection $N^1_i\cap N^2_i$
is a toric subvariety which is isomorphic to the product of the
toric degenerations of polygon spaces with side lengths $(r_1,
...., r_{i-1})$ and $(r_i, ..., r_n)$ respectively. Analogously,
the intersection $N^2_i$ and $N^3_i$, if non-empty, is a toric
subvariety which is isomorphic to the product of the toric
degenerations of polygon spaces with side lengths $(r_1, ....,
r_{i})$ and $(r_{i+1}, ..., r_n)$ respectively. It is clear that
these divisors $\{ N^j_i\}$ for $j=1,2$, and $3$ and $2\le i\le
n-1$ generate $A^*(N_{\bf r})$. Although the additive structure on
homology of $M_{\bf r}$ can readily be understood if we work with
divisors $\{ Z_{ij}\}$, there are too many generators and linear
relations between these cycles. Using the real cycles represented
by connected components of the $Y_{i}$'s, the homology of $M_{\bf
r}$ can be understood with more ease, as we can see explicitly the
toric subvarieties that these cycles degenerate to in $N_{\bf
r}$.}
\end{rem}

\begin{question} Do the cycles of the (connected components of)
subspaces $Y_{i}$ for $2\le i\le n-1$, generate $H_*(M_{\bf r},
\Z)$?
\end{question}

In some examples that we were able to compute, this indeed was the
case. However, in general due to the fact that some of the $Y_i$'s
degenerate to the sum of two or three irreducible toric divisors,
the answer to this question might not be entirely obvious.

As follows from \cite[Proposition 2.2(ii)]{AlexBrion}, our family
is isomorphic to the trivial family over $\C^*=\C\setminus\{ 0\}$.
In the future we hope to establish an explicit trivialization, so
the degeneration of the aforementioned cycles will be transparent.

\section{Degeneration of Flaschka-Millson integrable systems}

In \cite{FlMill} the authors considered moduli spaces of weighted
configurations on $\CP^m$. Given an $n$-tuple of positive integers
${\bf r}=(r_1, ..., r_n)$ their moduli space is, again, identified
with the GIT quotient of $(\CP^m)^n$ by the diagonal action of
${\rm SL}(m+1, \C)$, with the linearizing line bundle ${\mathcal
O}({\bf r})$ is given by the tensor product of pull-backs of the
bundles ${\mathcal O}(r_i)$ from the $i$-th factor. It has a
unique ${\rm SL}(m+1, \C)$ linearization. In order for this moduli space
to be nonempty, we need to impose, in addition to the obvious
condition $n
> m+1$, the so-called strong triangle inequalities, one for each
$1\le i\le n$:
$$
r_i\le P,\ \ \ \hbox{where}\ \ \ P:=\frac{1}{m+1}\sum_{j=1}^n r_j.
$$
Note that when $m=2$, these are just the usual triangle
inequalities. On the symplectic side, this moduli space is
identified with the moduli space of solution to the matrix
equation:
$$
A_1+\cdots + A_n=P.{\rm Id},
$$
where each $A_i$ is a rank one $(m+1)\times (m+1)$ hermitian
symmetric matrix, with spectrum $\{ r_i, 0, 0, ..., 0\}$, and we
mod out by the diagonal action of the group ${\rm U}(m+1)$. Then
the Flaschka-Millson bending hamiltonians, by definition, are the
non-trivial eigenvalues of the partial sums
$$
\sum_{i=1}^k A_i.
$$
Using the Gelfand-MacPherson correspondence \cite{Kapranov}, the
above moduli space $M_{\bf r}$ can be identified with a GIT
quotient of the grassmannian ${\rm Gr}(m+1, n)$ with respect to
the torus action. The choice of the linearized line bundle is
analogous to the case $m=1$.

These hamiltonians also appear if we take the Gelfand-Tsetlin
integrable system on the grassmannian and consider its reduction.
More precisely, consider the Gelfand-Tsetlin triangle for this
case. The top line consists of $(n-m-1)$ zeroes and $(m+1)$ $P$'s.
(For the purposes of drawing the diagram we assumed that $m< n/2$,
but we don't hold this assumption in general):
$$
\begin{array}{ccccccccccccccccccc}
0 & & 0 & & \cdots & & 0 & & 0 & & P & & P & & \cdots & & P & & P \\
 & 0 & & 0 & \cdots & 0 & & 0 & & a_1 & & P & & P & \cdots & P & & P & \\
 & & 0 & & \cdots & & 0 & & a_2 & & b_2 & & P & & \cdots & & P & & \\
 & & & \cdots & & & & & \cdots & & & & & & \cdots & & & & \\
& & & \cdots & & 0 & & a_{n-2m} & & \cdots & & p_{n-2m} & & P & & & & & \\
& & & & & \cdots & & & & \cdots & & & & & & & & & \\
& & & & 0 & & a_{m+1} & & & \cdots & & q_{m+1} & & & & & & & \\
& & & & & \cdots & & & \cdots & & \cdots & & & & & & & & \\
& & & & & & o_3 & & p_3 & & q_3 & & & & & & & & \\
& & & & & & & p_2 & & q_2 & & & & & & & & & \\
& & & & & & & & q_1 & & & & & & & & & & \\
\end{array}
$$
There are a total of $(m+1)\times (n-m-1)$ indeterminates $a_1$,
..., $q_{m+1}$ in the diagram. The reduction at the level ${\bf
r}=(r_1, ..., r_n)$ means that there are the following condition
imposed, analogous to Section 3, where we considered the case
$m=1$:
$$
q_1=r_1, \ \ p_2+q_2=r_1+r_2,\ \ ..., \ \ a_{m+1}+\cdots +
q_{m+1}=r_1+\cdots + r_{m+1}, \ ...\ \hbox{etc.}
$$
meaning that the sum of all elements in row $k$ from the bottom
equals $\sum_{i=1}^k r_i$. The Flaschka-Millson hamiltonians in
this presentation are the differences of adjacent elements in each
row:
$$
b_2-a_2, \ \ b_3-a_3, \ \ c_3-b_3, \ \ , ...., \ \ q_2-p_2.\
$$
There are precisely
$$
(m+1)(n-m-1)-(n-1)=mn-2m-m^2
$$
of them, which is exactly the complex dimension of the moduli
space. The polygon defined by the natural inequalities on these
hamiltonians is, as before, denoted by $\Pi_{\bf r}$.

Proceeding in a completely analogous way to Section 3, we can see
that there exists a flat degeneration of the moduli space $M_{\bf
r}$ to a polarized toric variety with the moment polytope
$\Pi_{\bf r}$.

\begin{rem} {\em If we denote by $P-{\bf r}$ the $n$-tuple of weights $(P-r_1, P-r_2, ..., P-r_n)$, then
Howard-Millson \cite{HowMill} showed that the GIT quotients
$$((\CP^m)^n//{\rm SL}(m+1, \C))_{{\mathcal O}({\bf r})}\ \ \ {\rm and} \ \ \
((\CP^{n-m-2})^n//{\rm SL}(n-m-1, \C))_{{\mathcal O}(P-{\bf r})}$$
are isomorphic. Applying the Gelfand-MacPherson correspondence,
this duality translates to the isomorphism of GIT quotients
$$({\rm Gr}(m+1, n)//H_1)_{{\mathcal O}(1), {\bf r}} \ \ \ {\rm and} \ \ \
({\rm Gr}(n-m-1, n)//H_1)_{{\mathcal O}(1),{P-{\bf r}}}\ \ , $$ which
they established by using the complex Hodge * operator. This
duality clearly induces an isomorphism between polytopes
$\Pi_{\bf r}$ and $\Pi_{P-{\bf r}}$ and the toric varieties
$N_{\bf r}$ and $N_{P-{\bf r}}$.}
\end{rem}

\section{Real loci}

First, let us have a brief discussion on real loci in general
weight varieties. Let $\tau$ be an involution of $G$ (group
automorphism of order $2$, and we treat $G$ as a real Lie group
with the Iwasawa decomposition $G=KAU$) satisfying the following
conditions:
\begin{itemize}
\item $\tau$ commutes with the Cartan involution $\theta$
\item $\tau$ maps $H$ to $H$
\item $\tau$ is anti-holomorphic
\item if $G^\tau$, $K^\tau$, and $A^\tau$ stand for the subgroups of $G$, $K$, and $A$ respectively fixed by
$\tau$ and $U^\tau$ is the subset of elements in $U$ fixed by
$\tau$, then $G^\tau=K^\tau A^\tau U^\tau$. (Note that we do not
require that $\tau$ maps $U$ to $U$.)
\end{itemize}
In fact, given a conjugacy class of real forms of $G$, one can
always find an involution $\tau$ as above so this conjugacy class
is represented by $G^\tau$. Alternatively, given a Satake diagram
$\Sigma$ \cite{Helga}, we can always find an involution $\tau$
such that $\Sigma$ corresponds to the real form $G^\tau$.

Now let us take a subset of simple roots $J$ such that it contains
all the black simple roots in the Satake diagram $\Sigma$ and
perhaps some white roots. If there are white roots connected by an
arrow, then we either include both of them to $J$ or neither. Let
$P_J$ be the corresponding parabolic subgroup containing $B$ (the
minimal among all parabolics containing $B$ and all the root
vectors $E_{-\alpha}$, $\alpha\in J$). The following is
straightforward:

\begin{lemma}The parabolic subgroup $P_J$ is $\tau$-stable.
\end{lemma}

This immediately implies that there is an induced involution, also
denoted by $\tau$, on the flag manifold $X=G/P$, and its fixed
point set is the real flag manifold $X^\tau\simeq G^\tau/P^\tau$.
In fact there exists an integral symplectic structure
\cite{FothReal} $\omega$ (corresponding to a polarizing line
bundle $L_\lambda$) such that $\tau$ is anti-symplectic with
respect to $\omega$ and $X^\tau$ is a lagrangian submanifold in
$X$.

In what follows we denote by  $T\subset H$ the connected component
of the maximal torus containing the identity element. The
involution $\tau$ maps $T$ to $T$ and let us denote by $T^+$ the
connected subgroup of $T$ on which $\tau$ acts identically and by
$T^-$ - where it acts by inverting all elements, so that
$T=T^+\times T^-$. Let us also denote by $Q$ a (not necessarily
connected) subgroup of $T$ consisting of {\it all} elements on
which $\tau$ acts by $\tau(t)=t^{-1}$. If the dimension of $T^+$
is denoted by $a$, then
$$
Q\simeq (\Z/2)^a\times T^-,
$$
and has $2^a$ connected components. Let us also denote by $T^\tau$
the subgroup of $T$ fixed by $\tau$. If the dimension of $T^{-}$
is $b$, then
$$
T^\tau\simeq T^+\times (\Z/2)^b.
$$

Note that the involution $\tau$ on $X$ satisfies the following
property with respect to $K$-action:
$$
\tau(k.x)=\tau(k)\tau(x).
$$
We can choose the moment map $\Phi$ for the $T$-action in such a
way that $\Phi(\tau(x))=-\tau(\Phi(x))$, where in the left hand
side the involution on $\ft^*$ is the one induced from $\tau$ on
$T$, and also denoted by $\tau$.  According to a theorem of
O'Shea-Sjamaar \cite{OSS} which generalizes Duistermaat's
\cite{Duister}, the image of $X^\tau$ under the moment map $\Phi$
for the $T$-action is the same as the intersection of $X$ with the
annihilator $(\ft^-)^\perp\subset \ft^*$ of $\ft^-$. Let us choose
a regular value $\mu$ of $\Phi$ such that $\mu\in(\ft^-)^\perp$.
Then the fixed point set of $\tau$ in $\Phi^{-1}(\mu)$, which is
the same as $\Phi^{-1}(\mu)\cap X^\tau$ is non-empty and
$T^\tau$-stable. If we reduce at $\mu$, and denote the induced
involution on $X//T$ also by $\tau$, then the natural map
$$
\psi:\ X^\tau//T^\tau:=(\Phi^{-1}(\mu)\cap X^\tau)/T^\tau \to
(X//T)^\tau
$$
has the following properties. The map $\psi$ is surjective onto a
(finite number of) connected component of $(X//T)^\tau$. The map
$\psi$ is a finite map, and is injective if $T$ acts freely on
$\Phi^{-1}(\mu)$. Replacing the involution $\tau$ acting on $X$ by
$s\tau$, where $s\in Q$ we can get all other connected components
of $(X//T)^\tau$ be in the image of maps analogous to $\psi$. The
induced involution on $X//T$ will still be the same, yet if $s$
acts non-trivially on $X$ and belongs to a different connected
component of $Q$ than the identity, then the $T^\tau$-orbits on
$X^{s\tau}$ are actually disjoint from those on $X^\tau$ and if
$X^{s\tau}$ is non-empty, then $X^{s\tau}\cap
\Phi^{-1}(\mu)/T^\tau$ will map onto different connected
components of $(X//T)^\tau$. To get all the connected components
of $(X//T)^\tau$, it is enough to choose one $s$ from each
connected component of $Q$, see details in \cite{FothReal}.

When the Satake diagram $\Sigma$ does not contain any black roots
(when $G^\tau$ is so-called quasi-split), then $B$ itself is
$\tau$-stable and the action of $\tau$ descends to the full flag
manifold $G/B$. And when $G^\tau$ is actually a split real form
(no arrows in $\Sigma$), then each $P\supset B$ is $\tau$-stable
and, moreover, $T=Q=T^{-}$ and $T^\tau$ is finite and isomorphic
to $(\Z/2)^r$, where $r={\rm rank}(G)$.

We will postpone a more detailed treatment of real loci in general
weight varieties and their degenerations to a future paper, and
now will concentrate on the case when $G={\rm SL}(n, \C)$, $\tau$
is the standard complex conjugation and $X$ is a complex
grassmannian. First, let us consider the case of spatial polygons.
Then the real quotient $X^{\tau}//T^\tau$ is the moduli space of
planar polygons (more precisely, its quotient by the mirror
reflections), which we denote by $M^{(2)}_{\bf r}$. If we, in
addition, fix admissible non-zero values of diagonals, then a
polygon can have only finitely many shapes. Their number
generically equals $2^{n-3}$, because we can do a 180${}^o$
``bending'' about each of the diagonals. Thus we have shown that
the map $M^{(2)}_{\bf r}\to \Pi_{\bf r}$ is surjective, and
generically $2^{n-3}$-to-one. However, this map is not globally
finite, and only becomes such at the special fiber $N_{\bf r}$. It
is easy to see that the involution $\tau$ on general fiber $M_{\bf
r}$ extends to a complex conjugate involution on the special
fiber, where its fixed point set is a lagrangian locus
$N^{(2)}_{\bf r}$ which maps surjectively and finitely onto
$\Pi_{\bf r}$. This statement has a natural generalization to the
moduli spaces $M_{\bf r}=(\CP^m)^n//{\rm SL}(m+1, \C)$ considered
in \cite{FlMill}:

\begin{proposition}
Let ${\bf r}=(r_1, ..., r_n)$ be an admissible $n$-tuple of
positive numbers. Let $\tau$ be a complex conjugate involution,
with fixed point set ${\rm SL}(n, \R)$ in ${\rm SL}(n, \C)$ and
${\rm Gr}_{\R}(m+1, n)$ in ${\rm Gr}_{\C}(m+1, n)$. Then $\tau$
extends to the special fiber $N_{\bf r}$ of the flat family and
the fixed point set $N^{\tau}_{\bf r}$ maps finitely and
surjectively under the moment map to the polytope $\Pi_{\bf r}$,
and the generic fiber has cardinality $2^{mn-2m-m^2}$.
\end{proposition}

\section*{Acknowledgments}
We thank Ben Howard and John Millson for showing us a preliminary
version of \cite{HowMill} and useful correspondence. We also thank
Michel Brion for a useful comment. The second author was partially
supported by an NSF grant.


\end{document}